\documentclass[12pt]{article}
\usepackage{epsf}
\usepackage{amsmath}
\usepackage{amssymb}

\newtheorem{preExample}{\bf Example}
\newenvironment{Example}{\begin{preExample}\hspace{-0.5em}. }{\end{preExample}}

\newtheorem{preLemma}{\bf Lemma}
\newenvironment{Lemma}{\begin{preLemma}\hspace{-0.5em}. }{\end{preLemma}}
\newcommand{\lem}[1]{\begin{Lemma}{#1}\end{Lemma}}

\newtheorem{preTHeorem}{\bf Theorem}
\newenvironment{THeorem}{\begin{preTHeorem}\hspace{-0.5em}. }{\end{preTHeorem}}
\newcommand{\thm}[1]{\begin{THeorem}{#1}\end{THeorem}}

\newtheorem{preCorollary}{\bf Corollary}
\newenvironment{Corollary}{\begin{preCorollary}\hspace{-0.5em}. }{\end{preCorollary}}
\newcommand{\corol}[1]{\begin{Corollary}{#1}\end{Corollary}}

\newtheorem{preDEfinition}{\bf Definition}
\newenvironment{DEfinition}{\begin{preDEfinition}\hspace{-0.5em}. }{\end{preDEfinition}}

\newenvironment{proof}{{\bf proof.}}{$\blacksquare$}
\newcommand{\prf}[1]{\begin{proof}#1\end{proof}}

\title{
An interesting proof of the nonexistence continuous bijection between $\mathbb{R}^n$ and $\mathbb{R}^2$ for $n\neq 2$ }
\author{ Fereshteh Malek
\thanks{Faculty of Sience, K.N.Toosi  University of Technology.
Email: {\tt malek@kntu.ac.ir}. }
\and Hamed Daneshpajouh 
\and Hamid Reza Daneshpajouh
\thanks{Department of Mathematics and Computer Science, University of Tehran. 
Email: {\tt h.r.daneshpajouh@khayam.ut.ac.ir}. }
\and  Johannes Hahn
\thanks{Faculty of Mathematics and Natural Sciences, University of Rostock.} }

\newcommand{\norm}[1]{\lVert#1\rVert}
\begin{document}
\maketitle
\begin{abstract}
\end{abstract}
In this article it is shown that there is no continuous bijection from $\mathbb{R}^n$ onto $\mathbb{R}^2$ for $n\neq 2$
by an elementary method. This proof is based on showing that for any cardinal number $\beta\leq 2^{\aleph_0}$, there is a partition of $R^n$  ($n\geq 3$) into $\beta$ 
arcwise connected dense subsets. 
\section{Introduction}
 In 1877 Cantor discovered a bijection of  $\mathbb{R}$  onto  $\mathbb{R}^n$,  
for any  $n\in\mathbb{N}$. Cantor's map was discontinuous, but the discovery of
Peano curve in 1890,  showed that there existed continuous (although not injective)
maps of  $\mathbb{R}$  onto  $\mathbb{R}^n$. After that and before 1910, several mathematicians showed that there didn't exist 
a bicontinuous bijection (homeomorphism) from  $\mathbb{R}^m$ onto  $\mathbb{R}^n$, 
for the cases $m=2$ and $m=3$ and $n>m$. Finally in 1911, Brouwer showed that there didn't exist a homeomorphism
 between $\mathbb{R}^m$ and  $\mathbb{R}^n$ for $n\neq m$ (For a modern treatment see Munkres, James (1984), p.109 [2]). 
The present paper proves not existance of continous bijection from  $\mathbb{R}^n$ onto  $\mathbb{R}^2$ for $n\neq 2$ by an elementary method.
 
Marey Ellen Rudin showed [1] that for any countable
cardinal $\alpha>2$, we can not partition the plane into $\alpha$ arcwise connected dense subsets. 
In this paper we show that for any cardinal number $\beta\leq 2^{\aleph_0}$, 
there is a partition of    $\mathbb{R}^n$ ($n\geq 3$ ) into $\beta$ arcwise connected dense subsets, and then by using this
we show that there is no continuous bijection from $\mathbb{R}^n$ onto $\mathbb{R}^2$, for $n\neq 2$. 
 
\lem{
There is a partition of $\mathbb R^{+} $ into $2^{\aleph_0}$ dense subsets. 
}
\prf{ Consider the additive group ($\mathbb R$,+). 
The quotient group $\mathbb R/\mathbb Q$ has $2^{\aleph_0}$ elements which are  dense subsets of  $\mathbb R$. 
Intersect them with $\mathbb R^{+}$. 
}
\thm{
There is a partition of $\mathbb R^3$ into $2^{\aleph_0}$  arcwise connected dense subsets. 
}
\prf{ 
Let $\{ S_i \;\vert\; i\in I\}$ be a  partition of  $\mathbb R^{+}$ into $2^{\aleph_0}$ dense subsets. 
$I$ is just an index set, so we may suppose that $I=(0\;  1)$.  Define 
 $L_i=\{(t, it , 0)\vert t> 0\} $ and  $M=\cup_{i\in I} L_i$ and
let $A_i$ be the union of all spheres  with center at the origin and  radius from  $S_i$, 
i.e 
$A_i=\{x\in\mathbb R^3\;\vert\; \norm{x}\in S_i \}$. 
Let $B_i=(A_i\backslash M)\cup L_i$. If $S$ is a sphere centered at the origin, then $S\backslash M $ is a sphere with a small arc removed, therefore
$A_i \backslash M $ is the union of some arcwise connected punctured spheres, 
open half-line $L_i$ pastes these punctured spheres together, so $B_i$ is arcwise connected.  
It is obvious that $\{B_i \vert i\in I\}$ is a partition of $\mathbb R^3$ with size $2^{\aleph_0}$. 
Since $S_i$ is dense in   $\mathbb R^{+}$, $A_i$ and consequently 
$B_i$ are dense in  $\mathbb R^3$. 
}
\corol{
There is a partition of $\mathbb R^n$ into $2^{\aleph_0}$  arcwise connected dense subsets for $n\geq 3$. 
}
\prf{ It is enough to set $B_i^{(n)}=B_i\times\mathbb R^{n-3} $,  
in which $B_i$ is as above.  
$\{B_i^{(n)}\vert i\in I\}$ is a partition of $\mathbb R^n$ which satisfies the claim.  
}

Note that the union of any number of the sets $B_i^{(n)}$ is an arcwise 
connected dense subset of $\mathbb R^n$, hence 

\corol{
For any cardinal number $\beta\leq2^{\aleph_0}$, there is a partition of $\mathbb R^n$ ($n\geq 3$) into $\beta$ 
arcwise connected dense subsets. 
}
\thm{
For any countable cardinal $\alpha > 2$ we can not partition the plane
 into $\alpha$ arcwise connected dense subset. 
}
\prf{
This statement is proved in [1]
}
\lem{
Let $X$,$Y$ be  metric spaces and $T:X\to Y$ be a continuous map

$(a)$ If $A$ is dense in $X$ and $T$ is onto, then $T(A)$ is dense in $Y$.
 
$(b)$ If $B\subset X$  is arcwise connected, then $T(B)$ is also arcwise connected.
}
\thm{
There is no continuous bijection from $\mathbb{R}$ onto $\mathbb{R}^m$ for $m\neq 1$.
}
\prf{
Suppose the contrary, let $g:\mathbb{R}\to\mathbb{R}^m$ be a continuous bijective map, we put $B_n=[-n,n]$ and so we have $\mathbb{R}^m=g(\cup_{n=1}^{\infty} B_n)=\cup_{n=1}^{\infty} g(B_n)$. $\mathbb{R}^m$ is not first category so at least one of the $g(B_n)$, forexample  $g(B_k)$ has a nonempty interior in $\mathbb{R}^m$, suppose $B(x,r)\subset g(B_k)$. Now we consider $f$ as a restriction of $g$ to $B_k$, since $B_k$  is compact so   $f:B_k\to g(B_k)$ is a homeomorphism and then  $B(x,r)$ is homeomorphic with an interval in $\mathbb{R}$ and that is a contradiction, because if we remove 3 points from $B(x,r)$ it remains connected but this is not the case for the intervals in  $\mathbb{R}$
}
\thm{
There is no continuous bijection from $\mathbb{R}^n$ onto  $\mathbb{R}^2$ for $n\neq 2$
}
\prf{
Suppose the contrary: 

(a) If $n>2$ then according to corollary $2$ and lemma $2$ we can partition $\mathbb{R}^2$ into 3 arcwise connected dense subsets and this contradicts theorem $2$. 

(b)If $n=1$  then this contradicts theorem $3$. 
}

\bigskip\noindent {\bf Acknowledgments.}
The authors are grateful to the professor Nicolas Hadjisavvas, for his valuable 
advices and comments.

\end{document}